# Demonstrating Efficiency Gains from Installing Truck Turntables at Crushers

M Brazil[1], P A Grossman[2], J H Rubinstein[3] and D A Thomas[4]

## ABSTRACT

The installation of truck turntables at above-ground crushers has the potential to yield gains in efficiency and productivity by eliminating the need for the trucks to turn around and reverse up to the crusher to dump their loads. The benefits include savings in time, fuel consumption and tyre wear. In addition, having a smaller area near the crushers dedicated to manoeuvring the trucks can be a significant benefit at some sites. Safety is also enhanced using turntables.

This paper describes a project with the Australian Turntable Company (ATC) to develop software that generates near-optimal layouts for the paths taken by the trucks in the vicinity of the crushers when turntables are present or absent. The geometry of the paths is constrained by the turning radius of the trucks and by the boundary of the manoeuvring area. The software quantifies the benefit of installing turntables by generating a design for paths with turntables and comparing it with the layout in current or planned use by the mining company. Through the use of this software tool, a mining company can readily assess the benefits of using turntables at their site.

## INTRODUCTION

Over the past ten to fifteen years there has been an increasing acknowledgement by the mining industry of the key role that mathematical optimisation can play in strategic mine planning and design. This has been demonstrated across a range of specific design tasks, from the original open pit optimal design framework pioneered by Lerchs and Grossmann (1965), to more recent advances in underground mine design including stope optimisation, the optimal design of underground access, and optimal scheduling of operations; see, for example, Alford, Brazil and Lee (2007) and Brazil *et al* (2009). The use of optimisation in mine planning has been shown not only to increase efficiency and reduce costs, but it can also significantly improve the overall net present value (NPV) of a mine.

One area of mine planning that has received very little attention to date is the design of truck transport routes from the mine to the above-ground crushers and waste dumps. The area in the vicinity of the crushers and waste dumps is often a source of traffic congestion, truck delays, and wear and tear on the trucks (due to the need for manoeuvring). This area also requires a large amount of cleared level space to accommodate truck movements, which may be a significant issue in protected environments or mountainous terrain. All of these factors can impact strongly on the NPV of the mine.

The installation of truck turntables at crushers and waste dumps has the potential to improve efficiency and productivity by eliminating the need for the trucks to turn around and reverse to dump their loads. While the idea and the use of turntables have been around since the 1800s, the use of truck turntables in industrial settings is relatively new. Nevertheless, it has already been shown to provide great benefits on building sites and in warehouses (for example, turntables for delivery trucks are now being extensively used by the Coles supermarket chain). The potential benefits of turntables on mining sites include savings in time, fuel consumption and tyre wear. In addition, since a smaller area is required near the crushers for manoeuvring the trucks, the need for cleared level space is substantially reduced. Finally, turntables can enhance safety in operations.

In order to rigorously assess the value of the truck turntable as an element of strategic mine designs it is important to be able to model and optimise truck movements in the vicinity of a crusher or waste dump. The Australian Turntable Company (ATC) manufactures and installs turntables for a variety of applications including turntables for ore trucks. In 2013, ATC approached the authors with a request for a software tool that would quantify the benefits of installing turntables at mine sites. The tool would be capable of rapidly generating, evaluating and comparing designs for the haulage paths taken by the trucks. The technical staff at ATC currently performs these tasks using CAD tools, but this is a time-consuming process and there is no guarantee that the designs are optimal. With the use of an automated tool, not only could ATC staff develop optimal designs efficiently but they could also quickly and easily demonstrate to mining companies the benefits of installing turntables at their mine sites.

This paper contains an outline of the mathematical framework used to solve this problem, a description of the

1. Associate Professor and Reader, Department of Electrical and Electronic Engineering, The University of Melbourne, Parkville Vic 3010. Email: brazil@unimelb.edu.au

2. Senior Research Fellow, Department of Mechanical Engineering, The University of Melbourne, Parkville Vic 3010. Email: peterag@unimelb.edu.au

3. MAusIMM, Professor, Department of Mathematics and Statistics, The University of Melbourne, Parkville Vic 3010. Email: rubin@ms.unimelb.edu.au

4. Professor, Head of Department and Associate Dean Research, Department of Mechanical Engineering, The University of Melbourne, Parkville Vic 3010. Email: doreen.thomas@unimelb.edu.au

Truck Path Optimisation Tool software that was developed based on this framework, a demonstration of its application using real data from a mine, and a discussion of planned future enhancements to the tool.

## DESCRIPTION OF THE TASK

The task can be best described in the context of a specific example. An image of the crusher area for a mine is shown in Figure 1. (For reasons of commercial confidentiality the mine is not identified here.) The crusher is indicated and the haulage paths in current use are clearly visible. Tipping can be performed from either the front or the rear of the crusher. Trucks enter from and leave towards two mining areas, one in the direction of the lower left and one in the direction of the right of the image. It is assumed that each truck returns to the same mining area as it departed from. Thus there are four haulage paths in this example: from and to each of the two mining areas via each of the two dump points.

The current practice in the absence of turntables is for the trucks to drive to a point near the crusher and then to stop, turn and reverse up to the crusher. The costs in time, fuel and tyre wear incurred by undertaking this manoeuvre are considerable. In addition, the manoeuvre can cause rocks to spill from the truck, resulting in further tyre wear for the trucks that follow.

The Truck Path Optimisation Tool requires a number of input parameters whose values depend on the type of trucks being used: kinematic parameters, fuel consumption rates, tyre wear rates and the turning circle radius. Parameters relating to the turntables are also needed. The parameters are listed in Table 1 along with the values they take in the present example. The values of the turntable parameters were supplied by ATC and the values of the other parameters were obtained from the mining company or from the published truck specifications.

**TABLE 1**
The input parameters for the Truck Path Optimisation Tool.

| | |
|---|---|
| Maximum forward truck speed | 10 k/h |
| Maximum reverse truck speed | 2.5 k/h |
| Truck acceleration | 0.5 m/s$^2$ |
| Truck deceleration | 1.8 m/s$^2$ |
| Tipping duration | 40 s |
| Turning circle radius | 28.4 m |
| Turntable maximum angular speed | 6 deg/s |
| Turntable angular acceleration/deceleration | 1.2 deg/s$^2$ |
| Turntable diameter | 15 m |
| Fuel consumption forward at maximum speed | 150 L/h |
| Fuel consumption reverse at maximum speed | 205 L/h |
| Fuel consumption accelerating forward unloaded | 361 L/h |
| Fuel consumption accelerating in reverse loaded | 395 L/h |
| Fuel consumption decelerating/idling | 53.7 L/h |
| Fuel consumption tipping | 211.7 L/h |
| Tyre wear loaded | 0.0231 mm/h |
| Tyre wear empty | 0.0119 mm/h |

The locations and bearings of the points where the trucks enter and leave the crusher area and the locations and exit bearings for the dump points are also required inputs to the program. In order to achieve the aim of keeping the program simple to use, it was decided that in the first version, the user should enter this information by clicking on the image.

Dumping can be performed at the front and the rear of the crusher. If turntables are in use then it is assumed that the turntable rotates before but not after tipping. Thus the truck can enter the turntable from any direction in a 180° range, but it always leaves the turntable by initially driving directly away from the crusher. The situation is illustrated in Figure 2. The diagram on the left shows a turntable in front of the crusher with a path that enters the turntable by the most direct route, which falls within the permitted 180° range. The diagram on the right shows a turntable at the rear of the crusher. In this case a direct entry path would fall outside the 180° range indicated and so instead the path enters the turntable from the left at right angles to the exit direction.

The Truck Path Optimisation Tool designs the haulage paths as minimum length curvature constrained paths. The theory of such paths is presented in the next section.

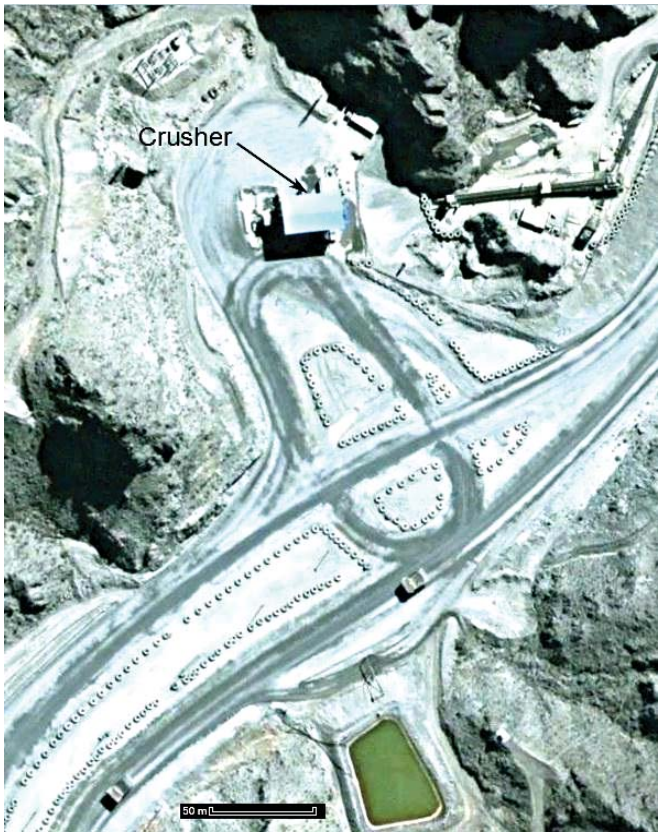

**FIG 1** – The crusher area.

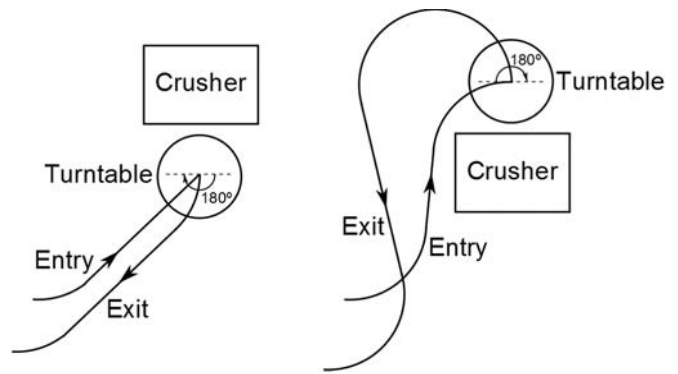

**FIG 2** – Truck movements for turntables in front of and behind the crusher.

## THE DESIGN OF OPTIMAL CURVATURE CONSTRAINED PATHS

The paths of the trucks are constrained by the minimum radius of their turning circle. In order to minimise the time taken, the fuel consumption and the tyre wear, the paths should be as short as possible subject to this constraint. Provided the truck is free to go anywhere in the area near the crusher, the optimal path for a truck between two given directed points is the shortest curvature-constrained path in the plane between the points. The problem of determining such paths was first investigated by Dubins (1957). Dubins showed that the optimal path is comprised entirely of circular arcs with minimum radius and straight line segments and that it must take one of the following forms:

- an arc followed by a line segment followed by an arc
- three arcs with alternating orientations
- a degenerate form of one of the above two forms in which one or more of the components have zero length.

Denoting a left turning arc, a right turning arc and a line segment by L, R and S respectively, and disregarding the degenerate forms, it follows that there are up to six potential minimum length paths between two given directed points: LSL, LSR, RSL, RSR, LRL and RLR. These paths are called Dubins paths. Some examples of Dubins paths are shown in Figure 3. It should be noted that not all of these Dubins paths need exist between a given pair of directed points. The optimal path can be obtained by simply finding all the Dubins paths between the given points and identifying the shortest one.

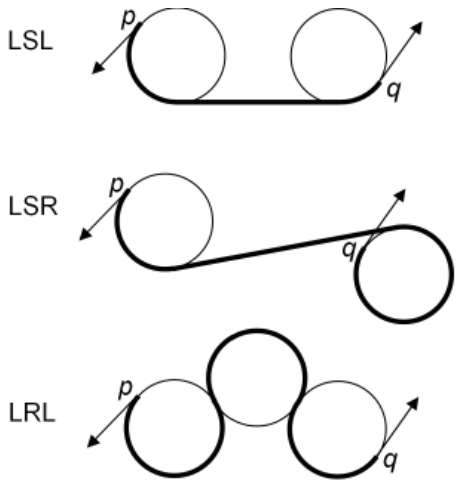

**FIG 3** – Examples of Dubins paths joining the directed points $p$ and $q$.

The authors have previously applied the theory of Dubins paths to the generation of optimal designs for the declines in underground mines (Brazil *et al*, 2009).

The Truck Path Optimisation Tool determines the shortest Dubins path of the form LSL, LSR, RSL or RSR between pairs of directed points specified by the user, except for the section of the path that enters the dump point which is handled as a special case. The Dubins paths with the forms LRL and RLR are not considered, since it is highly unlikely that a path with one of these forms would be optimal in the present context.

If there is a turntable at the dump point, then the section of the path that enters the dump point is handled differently because the entry bearing at the turntable is not specified uniquely but is merely restricted to a 180° range. First, the Truck Path Optimisation Tool generates the paths comprising an arc followed by a line segment (ie the paths of the form LS or RS). If at least one of these paths enters the turntable within the allowed 180° range then the shortest such path is chosen. If, on the other hand, neither of these paths enters the turntable within the allowed range then the Truck Path Optimisation Tool generates all of the LSL, LSR, RSL and RSR paths that enter the turntable at right angles to the specified exit direction from the turntable and the shortest of these paths is chosen.

If no turntable is present then a truck approaching a dump point must stop, turn and reverse up to the crusher. In order for the software tool to be able to generate and evaluate such paths, a 'reverse point' where the truck stops and then starts reversing must be inserted into the section of the path that precedes the dump point. The Truck Path Optimisation Tool creates a reverse point and determines its initial location and bearing using the theory of optimal curvature constrained paths with reversals developed by Reeds and Shepp (1990) and described below. The tool provides the facility for the user to change the location and bearing of the reverse point if required.

Reeds and Shepp (1990) investigated the problem of determining the shortest path in the plane between two directed points, subject to an upper bound on the curvature and with any number of cusps (corresponding to reverse points) permitted. They obtained a list of forms of such paths as sequences of line segments and circular arcs, analogous to the forms of the Dubins paths. For the present problem, only forms containing exactly one cusp need be considered, and then only those paths in which the cusp occurs just before the end of the path. Denoting the cusp by a vertical bar, the path forms that meet these criteria are LSL | R, LSR | L, RSL | R, RSR | L, LR | L and RL | R, where the orientation of the final (reverse) arc is relative to the direction the vehicle is facing and not the direction of travel. In addition, the second arc in each of the first four forms must turn through an angle of 90°. Examples of these paths are shown in Figure 4. As in the case of Dubins paths, the two path forms consisting of just three arcs are unlikely to arise in the present problem and can be disregarded. Given the start and end points, the optimal path can be found by determining the shortest path of the remaining four path forms.

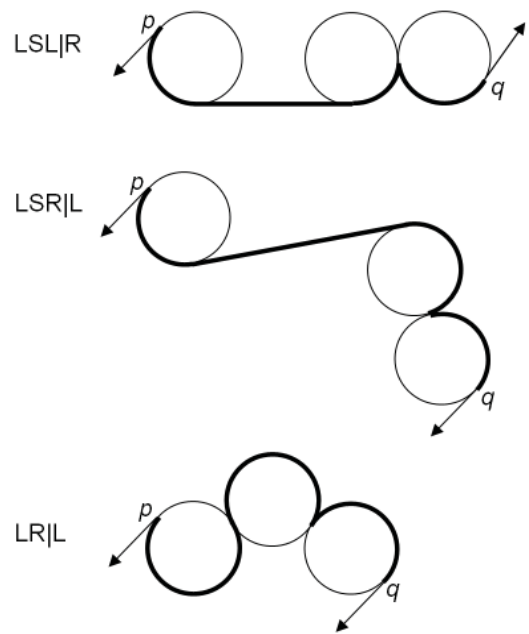

**FIG 4** – Examples of optimal paths with cusps following Reeds and Shepp (1990).

## THE TRUCK PATH OPTIMISATION TOOL

The Truck Path Optimisation Tool is a stand-alone program that runs under Windows. It was written in C++ and developed in Microsoft Visual Studio. The program is designed to be straightforward and intuitive to use.

After launching the program, the user loads an image of the dump area onto the screen. The image could be a map or diagram supplied by the mining company or a satellite image obtained from a source such as Google Maps. The user sets the scale of the image by clicking on two points on the image that are a specified distance apart. In the example in Figure 1, this would be done by clicking on each end of the scale bar in the image.

Next, the user would typically open a panel in which the parameters are displayed and change any values if required.

The user specifies the locations and bearings of the entry and exit points and the dump points by clicking on the image. The bearing of a dump point is the direction in which the truck leaves the dump point. Entry and exit points are input in pairs because it is assumed that each truck returns to the location from which it came. Any number of entry/exit point pairs and dump points may be input. The program generates a haulage path with a turntable and a haulage path without a turntable for every combination of an entry/exit pair and a dump point. The total time taken, the fuel consumption and the tyre wear for each path, and the savings obtained in time, fuel and tyre wear by installing turntables, are calculated and displayed.

Figures 5 and 6 show the main screen of the Truck Path Optimisation Tool with two entry/exit pairs and two dump points, one at the front of the crusher and one at the rear. The program has generated the four haulage paths shown on the image, one path for each combination of an entry/exit pair and a dump point, under the assumption that turntables are present (Figure 5) or absent (Figure 6). The time, fuel consumption and tyre wear for all the paths, with and without turntables, are tabulated at the bottom.

An enlarged image of the haulage paths with turntables is shown in Figure 7. In the terminology introduced earlier in context of Dubins paths, the forms of the entry/exit paths are LS/LSR, LS/LSL, LSL/LSR and RSL/RSL for the red, green, yellow and purple paths respectively.

The paths that are initially generated by the Truck Path Optimisation Tool sometimes pass through regions where trucks are unable to go. In order to handle this problem, the program provides the facility for the user to add waypoints to the paths. A waypoint is a directed point inserted into a section of a path in order to guide the route the path takes in that section. A waypoint is input by specifying the section of the path in which it is to occur and then specifying its location and bearing in the same manner as for entry/exit points and dump points. Figures 8 and 9 show the modified paths (with and without turntables respectively) that the program generates after waypoints have been inserted to guide the paths so that they remain outside any no-go regions.

The effect of installing turntables for the example shown in the figures can be read from the table at the bottom of the window. It reveals savings of between 15 and 85 seconds in time, between 2.0 and 5.8 L of fuel, and between 0.0003 and 0.0007 mm of tyre wear per trip, depending on the path. The number of trips per eight-hour shift at this site was estimated to be 109, based on the target production rate and the tray capacity. At three shifts per day and operating 365 days per year, the average annual savings from installing turntables at the crusher are around 1666 hours of truck time, 480 000 L of fuel, and 57 mm of tyre wear (which corresponds to more than three-quarters of the total permitted wear on a tyre). In practice, the reduction in truck times can mean that fewer trucks are required, with consequent savings in capital and maintenance costs.

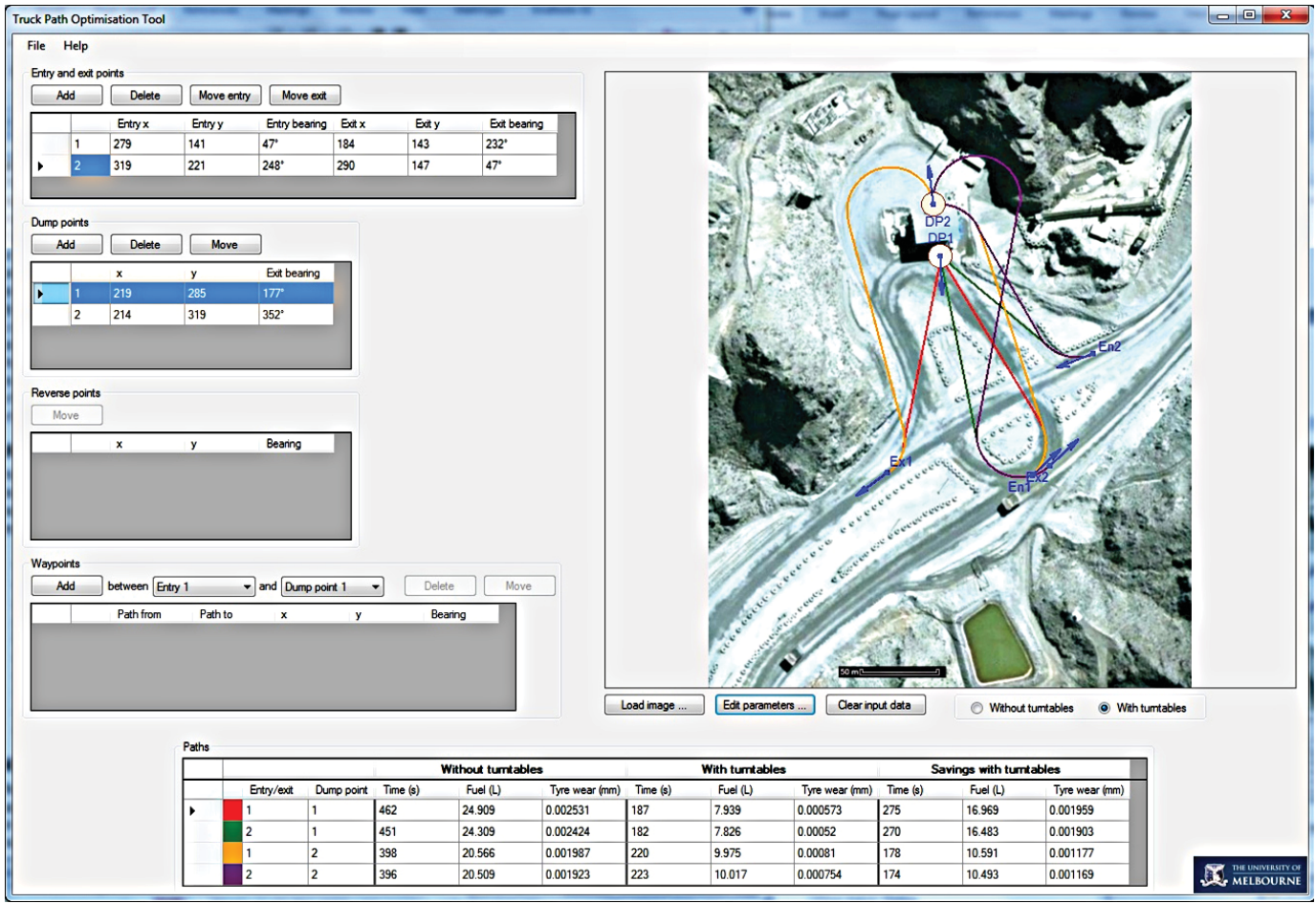

**FIG 5** – The user interface of the Truck Path Optimisation Tool showing path designs with turntables.

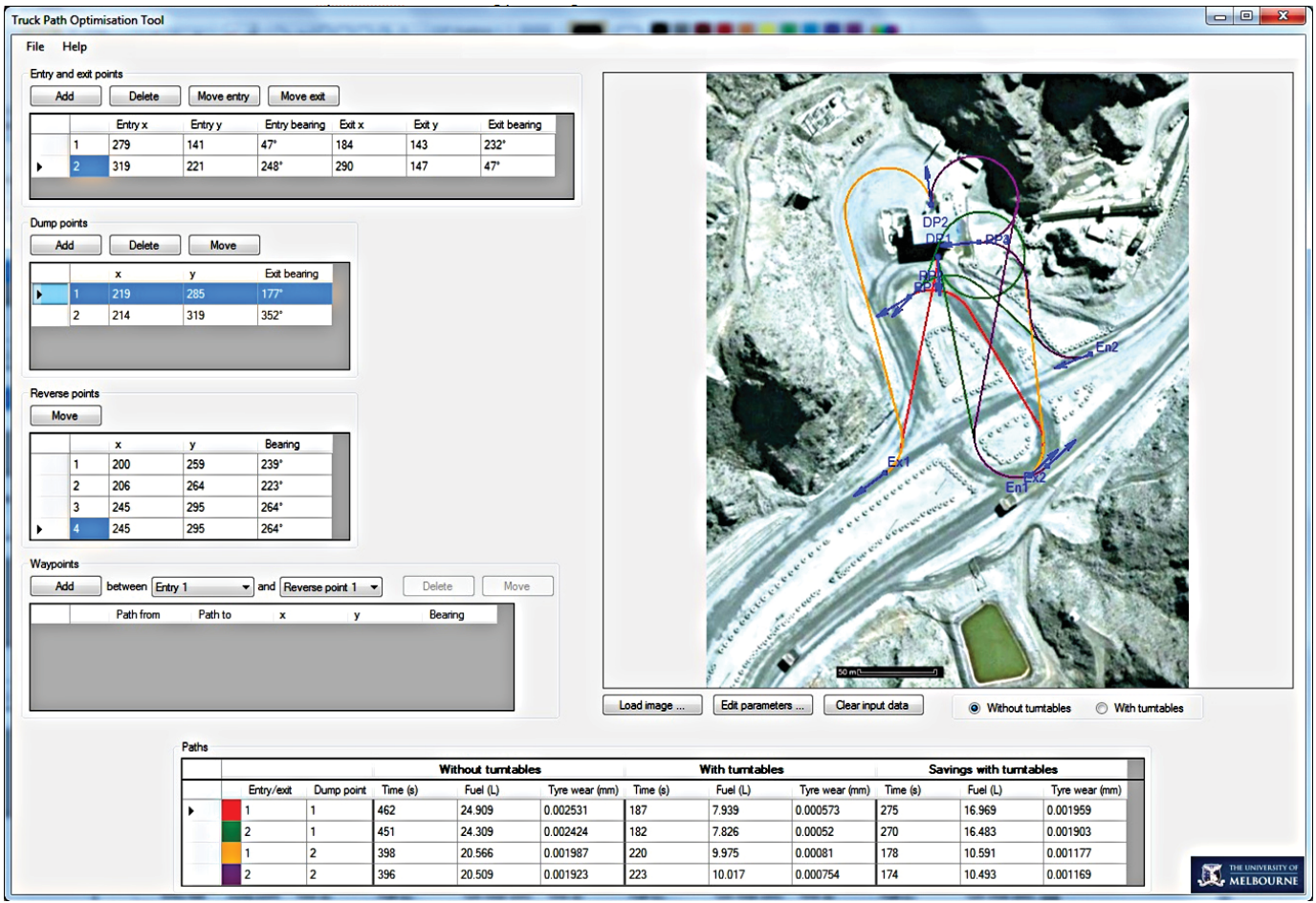

**FIG 6** – The user interface showing path designs without turntables.

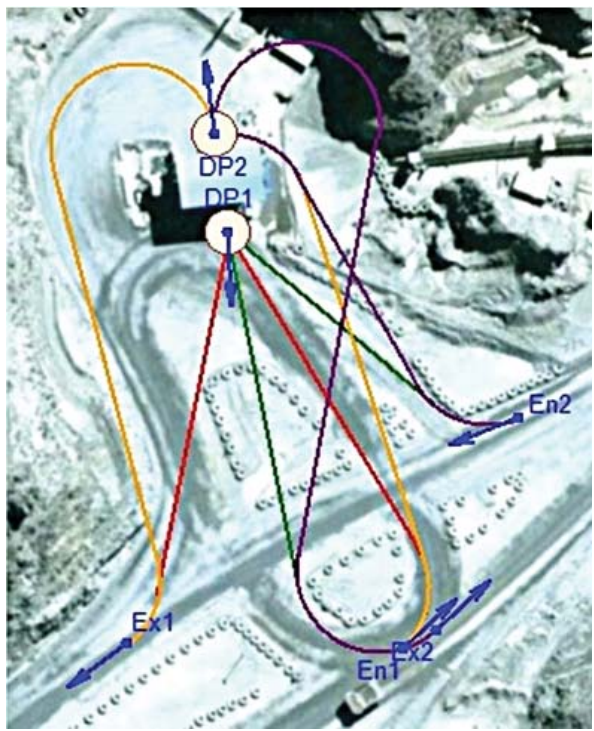

**FIG 7** – Detail from Figure 5 showing the haulage paths.

## PLANNED ENHANCEMENTS TO THE PROGRAM

The main deliverable of the project was an initial version of the Truck Path Optimisation Tool program with sufficient functionality to be useful to ATC. However, there are many ways in which the program could be developed further. If the opportunity becomes available, some or all of the following enhancements to the program are planned:

- The capability to calculate dollar costs for time, fuel and tyre wear. This capability would allow the 'bottom line' saving made by installing turntables to be calculated.

- A facility for the user to place barriers on the image to define no-go regions. The program would generate optimal paths that lie within the region outside the barriers.

- Modification of the paths to eliminate path crossings or, where that is not feasible, to ensure that paths cross at an acceptable angle. This would help to reduce delays caused by trucks having to give way.

- Enhancements to improve the usability of the program, including the capability of writing the input data and the generated outputs to files, and alternative methods of data entry such as keying in the point coordinates.

**FIG 8** – The user interface showing the path designs with turntables and with waypoints added to the paths.

**FIG 9** – The user interface showing the path designs without turntables and with waypoints added to the paths.

## CONCLUSIONS

A software tool, the Truck Path Optimisation Tool, has been developed to generate, evaluate and compare truck haulage paths at dump sites (crushers or waste dumps) where turntables are present or absent. The tool can be used to rapidly generate alternative designs and to demonstrate to clients the benefits of installing turntables. A case study was performed that provided an opportunity to test the tool using real data, while also demonstrating the benefits of installing turntables at the mine. Further work is planned to provide the tool with additional features including the capability of generating designs that take better account of crossing paths.

## ACKNOWLEDGEMENTS

The authors wish to thank Paul Chapman, Ben Chapman and Bassam El Aawar from the ATC for their help, support and enthusiasm throughout the project. ATC's approval to publish this material is also acknowledged. This project was financially supported by the Victorian Government Technology Voucher Program.